\documentclass[11pt]{article}

\usepackage[margin=1in]{geometry}
\usepackage[T1]{fontenc}
\usepackage{lmodern}
\usepackage{microtype}
\usepackage{amsmath,amssymb,amsthm,mathtools}
\usepackage{graphicx}
\usepackage{booktabs}
\usepackage{enumitem}
\usepackage{xurl}
\usepackage[hidelinks]{hyperref}
\usepackage[nameinlink,noabbrev]{cleveref}
\usepackage{listings}
\allowdisplaybreaks
\hypersetup{
  pdftitle={Theta-Eta Mellin Transforms, Shifted Reflections, and Rigorous Approximation Theory for a Variant Theta Construction},
  pdfauthor={Sourangshu Ghosh and Darrell Cox}
}

\newtheorem{theorem}{Theorem}[section]
\newtheorem{proposition}[theorem]{Proposition}
\newtheorem{corollary}[theorem]{Corollary}

\theoremstyle{definition}
\newtheorem{definition}[theorem]{Definition}
\newtheorem{remark}[theorem]{Remark}
\newtheorem{question}[theorem]{Research question}

\newcommand{\C}{\mathbb{C}}

\newcommand{\RePart}{\operatorname{Re}}

\newcommand{\etaD}{\eta}
\newcommand{\LamZ}{\Lambda_{\zeta}}
\newcommand{\LamE}{\Lambda_{\eta}}

\newcommand{\ThetaB}{\Theta^{B}}

\newcommand{\etaN}{\eta_N}

\newcommand{\Oof}{\mathcal{O}}

\title{Theta--Eta Mellin Transforms, Shifted Reflections, and\\Rigorous Approximation Theory for a Variant Theta Construction}
\author{Sourangshu Ghosh\\
\small Indian Institute of Science, Bengaluru, India
\and Darrell Cox\\
\small Department of Mathematics, Grayson County College, United States}
\date{}

\begin{document}
\maketitle

\begin{abstract}
We develop a rigorous analytic framework for investigating a proposed functional equation associated with a variant theta construction. We prove the classical Mellin identity connecting the Jacobi kernel $1-\theta_4(0\mid ix^2)$ to the completed Dirichlet eta function, and extend its Mellin transform meromorphically to the whole plane with a single simple pole at the origin. Mellin uniqueness identifies $x^{3/8}(1-\theta_4(0\mid ix^2))$ as the unique kernel in a natural weighted $L^1$ class whose Mellin transform is the shifted completed eta function. We distinguish this Gaussian kernel from the weighted exponential series $\Theta^B(s)=\sum_{n\geq1}ne^{-ns}$, whose Mellin transform is $\Gamma(s)\zeta(s-1)$ and obeys a different, explicitly derived reflection relation. We show that every affine reflection $s\mapsto c-s$ can be generated for a shifted completed eta or zeta function by choosing the shift appropriately; consequently, the occurrence of $\tfrac14-s$ alone cannot identify the original variant with a shifted eta function. The corrected analytic product defined here is an elementary meromorphic multiplier of the completed zeta function and therefore has a precise $1-s$ functional equation. We prove an explicit eta-tail bound, a sharper alternating-tail asymptotic, a uniformly convergent Euler-transformed eta series with a geometric error bound, a uniform zeta-approximation bound away from the line $\Re s=1$, and a combined truncation estimate when both the zeta and weighted-exponential factors are truncated. We also determine the principal Laurent data and the cancellation of the gamma poles at the trivial zeta zeros for the corrected product. Reproducible high-precision experiments verify the stated error scales. The exact identification of the original variant with a shifted completed eta function remains open; no claim of a new functional equation for that undefined legacy limit is made.
\end{abstract}


\section{Introduction and status of the problem}

The motivating study considers a finite complex-valued construction involving a gamma-type product, the weighted exponential series
\begin{equation}\label{eq:thetaB}
\ThetaB(t)=\sum_{n=1}^{\infty}ne^{-nt},\qquad \RePart(t)>0,
\end{equation}
and partial sums used as numerical approximations to the Riemann zeta function. The numerical study compares curves built from these ingredients and proposes replacing the usual reflected argument $1-s$ by $\tfrac14-s$ in an equation for the logarithm of the constructed quantity.

Three questions must be kept separate:
\begin{enumerate}[label=\textup{(Q\arabic*)},leftmargin=3em]
  \item What exact Mellin transform of a classical theta kernel gives the completed Dirichlet eta function?
  \item Under what definition does a reflection $s\mapsto\tfrac14-s$ follow rigorously?
  \item Which finite approximations converge to the stated limiting objects, and what quantitative errors accompany them?
\end{enumerate}

The classical identities answer these questions for explicitly defined functions, but they do not by themselves identify the distinct weighted exponential construction in \eqref{eq:thetaB} with the Jacobi theta kernel. The revised analysis therefore distinguishes the original ill-defined limiting proposal from a corrected analytic surrogate, and supplies an independent, reproducible numerical study with explicit error bounds.

\paragraph{Scope and contribution.}
The Jacobi theta--eta Mellin identity and the completed zeta functional equation are classical; detailed proofs are included to make normalizations and domains explicit. The additional results are structural consequences developed for this investigation: (i) meromorphic continuation and a finite-part integral for the theta-kernel Mellin transform; (ii) a distinct functional equation for the Mellin transform of $\Theta^B$; (iii) the full family of affine reflections generated by shifts; (iv) an exact reflection law for the corrected product; and (v) an explicit error bound when both factors in the finite product are truncated; (vi) an Euler-transformed eta expansion with a rigorous geometric tail bound; and (vii) a local singularity analysis of the corrected product, including its Laurent coefficients at $s=0,1$ and removable singularities at negative even integers; and (viii) a Mellin-transform uniqueness criterion that specifies the kernel required by the shifted eta identity. We also prove a sharper alternating-tail expansion explaining the observed $-1/2$ slope on the critical line. These statements do not identify the corrected product with the original legacy construction, whose nonreal Euler limit does not converge. Because the external zero list and raw samples underlying the old crossing-index figures were not supplied, those legacy regressions are not claimed to have been regenerated or validated. No priority or novelty claim is made for identities that follow from the classical zeta functional equation.

Throughout, $s=\sigma+it$, $x^s=e^{s\log x}$ for $x>0$, and all functional equations are meromorphic identities unless a domain is specified. We use the standard completed functions
\begin{equation}\label{eq:completed-definitions}
\LamZ(s)=\pi^{-s/2}\Gamma\!\left(\frac{s}{2}\right)\zeta(s),\qquad
\LamE(s)=\pi^{-s/2}\Gamma\!\left(\frac{s}{2}\right)\etaD(s).
\end{equation}

\section{Obstructions in the earlier finite formulation}

\subsection{The real-part-based Euler limit}

The ordinary Euler limit is
\begin{equation}\label{eq:euler}
\Gamma(s+1)=\lim_{N\to\infty}
\frac{N!(N+1)^s}{\prod_{k=1}^{N}(k+s)},\qquad \RePart(s)>-1.
\end{equation}
The earlier formulation instead replaces $s$ in the denominator by its real part. Define
\begin{equation}\label{eq:Pn}
P_N(s)=\frac{N!(N+1)^s}{\prod_{k=1}^{N}(k+\RePart(s))},\qquad s=\sigma+it.
\end{equation}
\begin{proposition}\label{prop:Pn}
If $\sigma>-1$ and $t\ne0$, then $P_N(s)$ has no limit as $N\to\infty$. If $t=0$, then $P_N(s)\to\Gamma(\sigma+1)$.
\end{proposition}
\begin{proof}
The gamma recurrence gives
\[
P_N(s)=\Gamma(\sigma+1)\frac{\Gamma(N+1)(N+1)^s}{\Gamma(N+\sigma+1)}.
\]
The gamma-ratio asymptotic yields
\[
P_N(s)=\Gamma(\sigma+1)N^{it}(1+o(1)).
\]
For $t=0$ this tends to $\Gamma(\sigma+1)$. For $t\ne0$, select integer sequences nearest to $\exp(2\pi k/|t|)$ and $\exp((2\pi k+\pi)/|t|)$. Along these sequences $N^{it}$ tends to $1$ and $-1$, respectively. The limit therefore does not exist.
\end{proof}

Thus the real-part-based limit cannot serve as a complex gamma factor. In this revision we use the ordinary gamma function in the explicitly defined corrected product in Section~\ref{sec:corrected-model}; we do not silently identify it with the nonconvergent limit \eqref{eq:Pn}.

\subsection{The weighted exponential series}

For $\RePart(t)>0$, summing a differentiated geometric series gives
\begin{equation}\label{eq:thetaBclosed}
\ThetaB(t)=\frac{e^{-t}}{(1-e^{-t})^2}.
\end{equation}
Its Mellin transform is distinct from the Jacobi theta transform.
\begin{proposition}\label{prop:thetaB-mellin}
For $\RePart(s)>2$,
\begin{equation}\label{eq:thetaB-mellin}
\int_0^\infty t^{s-1}\ThetaB(t)\,dt=\Gamma(s)\zeta(s-1).
\end{equation}
\end{proposition}
\begin{proof}
For $\RePart(s)>2$, absolute convergence permits interchange of the series and integral:
\[
\begin{aligned}
\int_0^\infty t^{s-1}\ThetaB(t)\,dt
&=\sum_{n=1}^{\infty}n\int_0^\infty t^{s-1}e^{-nt}\,dt\\
&=\Gamma(s)\sum_{n=1}^{\infty}n^{1-s}
=\Gamma(s)\zeta(s-1).
\end{aligned}
\]
\end{proof}

The identity \eqref{eq:thetaB-mellin} does not exclude a different transform identity for a larger construction. It establishes that the weighted series itself does not produce the eta Mellin transform under its natural Mellin integral. Any connection between the whole variant and the eta function requires an additional explicit identity.

\subsection{Periodicity, truncation, and the reflection of the weighted-series Mellin transform}
The closed form \eqref{eq:thetaBclosed} has the equivalent representation
\begin{equation}\label{eq:thetaB-hyperbolic}
\ThetaB(s)=\frac{1}{4\sinh^2(s/2)}.
\end{equation}
It follows that $\ThetaB$ is an even meromorphic function of period $2\pi i$, has no zeros, and has double poles precisely at $2\pi i\mathbb Z$. Thus, for fixed $\sigma>0$, its value along $s=\sigma+it$ depends only on $t$ modulo $2\pi$; this exact periodicity is a possible source of repeated patterns in plots involving the factor, and is unrelated by itself to a Jacobi theta--eta transformation. In particular, on each half-plane $\RePart(s)\geq\delta>0$ it is uniformly bounded in the imaginary part:
\begin{equation}\label{eq:thetaB-strip-bound}
|\ThetaB(s)|\leq \frac{e^{-\delta}}{(1-e^{-\delta})^2}.
\end{equation}
Indeed, with $q=e^{-s}$ and $|q|\leq e^{-\delta}$, one has $|1-q|\geq1-|q|\geq1-e^{-\delta}$.

For $N\geq1$, define the finite weighted sum $\Theta^B_N(s)=\sum_{n=1}^N ne^{-ns}$. The tail has an exact formula and a useful explicit bound.
\begin{proposition}[Weighted-exponential truncation bound]\label{prop:thetaB-tail}
Let $s=\sigma+it$ with $\sigma>0$ and $q=e^{-s}$. Then
\begin{align}
\ThetaB(s)-\Theta^B_N(s)
&=q^{N+1}\left(\frac{N+1}{1-q}+\frac{q}{(1-q)^2}\right),\label{eq:thetaB-tail-exact}\\
|\ThetaB(s)-\Theta^B_N(s)|
&\leq e^{-\sigma(N+1)}\left(\frac{N+1}{1-e^{-\sigma}}+\frac{e^{-\sigma}}{(1-e^{-\sigma})^2}\right).\label{eq:thetaB-tail-bound}
\end{align}
The convergence is locally uniform in $\RePart(s)>0$ and is exponentially fast on every closed half-plane $\RePart(s)\geq\delta>0$.
\end{proposition}
\begin{proof}
The geometric-series identity and its derivative give, for $|q|<1$,
\[
\sum_{n=N+1}^{\infty}nq^n=q^{N+1}\left(\frac{N+1}{1-q}+\frac{q}{(1-q)^2}\right),
\]
which proves \eqref{eq:thetaB-tail-exact}. The inequalities $|q|=e^{-\sigma}$ and $|1-q|\geq1-|q|$ imply \eqref{eq:thetaB-tail-bound}. On compact subsets of $\RePart(s)>0$, the right-hand side converges uniformly to zero; when $\sigma\geq\delta$, it is $O_\delta((N+1)e^{-\delta N})$.
\end{proof}

The Mellin transform in Proposition~\ref{prop:thetaB-mellin} also satisfies its own functional equation. Write $B(s)=\Gamma(s)\zeta(s-1)$ for its meromorphic continuation.
\begin{theorem}[Functional equation for the weighted-series Mellin transform]\label{thm:weighted-Mellin-FE}
Write $B(s)=\Gamma(s)\zeta(s-1)$ for the meromorphic continuation of the Mellin transform in Proposition~\ref{prop:thetaB-mellin}. Then, as a meromorphic identity,
\begin{equation}\label{eq:weighted-Mellin-FE}
\boxed{\sin\left(\frac{\pi s}{2}\right)B(s)
=(s-1)2^{s-2}\pi^{s-1}\zeta(2-s).}
\end{equation}
Equivalently,
\begin{equation}\label{eq:weighted-Mellin-reflection}
\boxed{B(s)=\frac{(s-1)2^{s-2}\pi^{s-1}}{\Gamma(3-s)\sin(\pi s/2)}B(3-s),}
\end{equation}
with the quotient interpreted meromorphically. Thus the reflection of the Mellin variable is $s\mapsto3-s$, whose fixed line is $\RePart(s)=3/2$. The zeta factor $\zeta(s-1)$ itself reflects to $\zeta(2-s)$; this one-unit difference in the argument must be retained when translating the zeta functional equation to the Mellin variable.
\end{theorem}
\begin{proof}
The standard zeta functional equation, applied at $s-1$, yields
\[
\zeta(s-1)=2^{s-1}\pi^{s-2}\sin\left(\frac{\pi(s-1)}2\right)\Gamma(2-s)\zeta(2-s).
\]
Multiplying by $\Gamma(s)$ and using
\[
\Gamma(s)\Gamma(2-s)=\frac{(1-s)\pi}{\sin(\pi s)},\qquad
\sin\left(\frac{\pi(s-1)}2\right)=-\cos\left(\frac{\pi s}{2}\right),
\]
then simplifying $\sin(\pi s)=2\sin(\pi s/2)\cos(\pi s/2)$ gives \,\eqref{eq:weighted-Mellin-FE} wherever the intermediate quotients are defined. Since $B(3-s)=\Gamma(3-s)\zeta(2-s)$, substitution yields \,\eqref{eq:weighted-Mellin-reflection}. Both cross-multiplied identities extend meromorphically to all $s$.
\end{proof}

\subsection{Ordinary partial zeta sums in the critical strip}

Let $\zeta_N(s)=\sum_{n=1}^N n^{-s}$. These sums converge to $\zeta(s)$ only for $\RePart(s)>1$. For $0<\sigma=\RePart(s)<1$, Euler--Maclaurin summation gives
\begin{equation}\label{eq:zetaN-asymptotic}
\zeta_N(s)=\frac{N^{1-s}}{1-s}+\zeta(s)+\Oof_s(N^{-\sigma}).
\end{equation}
Since $|N^{1-s}/(1-s)|=N^{1-\sigma}/|1-s|\to\infty$, the sequence diverges, including when $s$ is a nontrivial zero of $\zeta$. The finite sums are Dirichlet polynomials in this region; they are not convergent approximations to the analytically continued zeta function.

\subsection{The logarithm used in the legacy routine}

For $z=re^{i\vartheta}\ne0$, the principal complex logarithm is
\begin{equation}\label{eq:log-correct}
\operatorname{Log} z=\log r+i\operatorname{Arg} z.
\end{equation}
The legacy C routine computes the two components as $\log(r)\cos\vartheta$ and $\log(r)\sin\vartheta$, which represents $(\log r)e^{i\vartheta}$ rather than \eqref{eq:log-correct}. A correct principal-branch implementation would store the real part as $\log r$ and the imaginary part as $\vartheta$; a continuous phase plot would additionally need branch unwrapping. At $z=0$, the logarithm is undefined. The numerical experiments below use real quantities $\log_{10}|\text{error}|$ and do not treat them as complex logarithms.

\section{The classical theta--eta Mellin identity}\label{sec:theta-eta}

The Dirichlet eta function is
\begin{equation}\label{eq:eta-def}
\etaD(s)=\sum_{n=1}^{\infty}\frac{(-1)^{n-1}}{n^s}
=(1-2^{1-s})\zeta(s),\qquad \RePart(s)>0,
\end{equation}
where the series converges in the stated half-plane and the right-hand side gives its analytic continuation. The singularity of the product at $s=1$ is removable for $\etaD$, which is entire.

Use the standard Jacobi theta convention
\begin{equation}\label{eq:theta4-def}
\theta_4(0\mid ix^2)=1+2\sum_{n=1}^{\infty}(-1)^n e^{-\pi n^2x^2},\qquad x>0,
\end{equation}
and set
\begin{equation}\label{eq:kernel}
K(x)=1-\theta_4(0\mid ix^2)=2\sum_{n=1}^{\infty}(-1)^{n-1}e^{-\pi n^2x^2}.
\end{equation}
\begin{theorem}[Theta--eta Mellin identity]\label{thm:theta-eta}
For every $s\in\C$ with $\RePart(s)>0$, the integral converges absolutely and
\begin{equation}\label{eq:theta-eta}
\boxed{\int_0^\infty x^{s-1}\bigl(1-\theta_4(0\mid ix^2)\bigr)\,dx
=\pi^{-s/2}\Gamma\!\left(\frac{s}{2}\right)\etaD(s)=\LamE(s).}
\end{equation}
\end{theorem}
\begin{proof}
For $\RePart(s)>1$, absolute convergence justifies termwise integration in \eqref{eq:kernel}. Since
\[
\int_0^\infty x^{s-1}e^{-\pi n^2x^2}\,dx
=\frac12\pi^{-s/2}n^{-s}\Gamma\!\left(\frac{s}{2}\right),
\]
we obtain
\[
\begin{aligned}
\int_0^\infty x^{s-1}K(x)\,dx
&=2\sum_{n=1}^{\infty}(-1)^{n-1}\int_0^\infty x^{s-1}e^{-\pi n^2x^2}\,dx\\
&=\pi^{-s/2}\Gamma\!\left(\frac{s}{2}\right)\etaD(s).
\end{aligned}
\]
To extend the equality, use the Jacobi transformation
\[
\theta_4(0\mid ix^2)=x^{-1}\theta_2(0\mid i/x^2).
\]
As $x\downarrow0$, $\theta_2(0\mid i/x^2)=\Oof(e^{-\pi/(4x^2)})$, so $K(x)=1+\Oof(x^{-1}e^{-\pi/(4x^2)})$. As $x\to\infty$, $K(x)=\Oof(e^{-\pi x^2})$. Thus the Mellin integral converges absolutely for $\RePart(s)>0$ and is locally uniformly dominated on compact subsets of that half-plane. It defines a holomorphic function there. The right-hand side is also holomorphic there, and the identity theorem completes the proof. This classical formula is recorded in DLMF, Eq.~20.10.3.
\end{proof}

\begin{theorem}[Meromorphic continuation of the theta-kernel Mellin transform]\label{thm:theta-Mellin-continuation}
Let $K(x)=1-\theta_4(0\mid ix^2)$ and let $\mathcal M K(s)=\int_0^\infty x^{s-1}K(x)\,dx$ initially for $\RePart(s)>0$. Then
\begin{equation}\label{eq:theta-Mellin-continuation}
\mathcal M K(s)=\frac1s+E(s),\qquad E\text{ entire on }\mathbb C,
\end{equation}
where
\begin{equation}\label{eq:E-entire}
E(s)=\int_0^1x^{s-1}(K(x)-1)\,dx+\int_1^\infty x^{s-1}K(x)\,dx.
\end{equation}
Consequently, $\Lambda_\eta(s)$ has exactly one pole, a simple pole at $s=0$ with residue $1$. Moreover,
\begin{equation}\label{eq:finitepart-theta}
\operatorname{FP}_{s=0}\mathcal M K(s)
=\int_0^1\frac{K(x)-1}{x}\,dx+\int_1^\infty\frac{K(x)}{x}\,dx
=\frac{\log\pi-\gamma-2\log2}{2},
\end{equation}
where $\gamma$ is Euler's constant. At positive integers, in particular,
\begin{equation}\label{eq:theta-kernel-moments}
\int_0^\infty K(x)\,dx=\log2,\qquad
\int_0^\infty xK(x)\,dx=\frac{\pi}{12},\qquad
\int_0^\infty x^2K(x)\,dx=\frac{3\zeta(3)}{8\pi}.
\end{equation}
\end{theorem}
\begin{proof}
The Jacobi transformation used in Theorem~\ref{thm:theta-eta} gives $K(x)-1=-x^{-1}\theta_2(0\mid i/x^2)=O(x^{-1}e^{-\pi/(4x^2)})$ as $x\downarrow0$. At infinity, $K(x)=O(e^{-\pi x^2})$. These estimates imply that both integrals in \eqref{eq:E-entire} converge locally uniformly for every complex $s$, so $E$ is entire by standard parameter-dependent integration. For $\RePart(s)>0$, splitting the Mellin integral at $1$ gives \eqref{eq:theta-Mellin-continuation}; since it equals $\Lambda_\eta(s)$ there, it supplies the meromorphic continuation and shows that the only pole is at $0$ with residue $1$.

For the constant term, use $\eta(s)=(1-2^{1-s})\zeta(s)$ and the standard values $\zeta(0)=-1/2$ and $\zeta'(0)=-\tfrac12\log(2\pi)$ to obtain $\eta(0)=1/2$ and $\eta'(0)=\tfrac12\log(\pi/2)$. The expansions $\Gamma(s/2)=2/s-\gamma+O(s)$ and $\pi^{-s/2}=1-(s/2)\log\pi+O(s^2)$ yield
\[
\Lambda_\eta(s)=\frac1s+\frac{\log\pi-\gamma-2\log2}{2}+O(s).
\]
Comparing constant terms in \eqref{eq:theta-Mellin-continuation} proves \eqref{eq:finitepart-theta}. Finally, substitute $s=1,2,3$ in the Mellin identity and use $\eta(1)=\log2$, $\eta(2)=\pi^2/12$, and $\eta(3)=3\zeta(3)/4$ to get \eqref{eq:theta-kernel-moments}.
\end{proof}

\begin{theorem}[Mellin inversion for the Jacobi kernel]\label{thm:theta-Mellin-inversion}
For every $c>0$ and $x>0$,
\begin{equation}\label{eq:theta-Mellin-inversion}
K(x)=\frac{1}{2\pi i}\int_{c-i\infty}^{c+i\infty}\Lambda_\eta(s)x^{-s}\,ds,
\end{equation}
where the integral is absolutely convergent.
\end{theorem}
\begin{proof}
Set $f_c(y)=e^{cy}K(e^y)$. The small- and large-$x$ estimates used above show that $f_c$ and its first two derivatives are integrable on $\mathbb R$ and vanish at infinity. With the Fourier-transform convention $\widehat f(t)=\int_{\mathbb R}f(y)e^{ity}\,dy$, the substitution $x=e^y$ gives $\widehat f_c(t)=\mathcal M K(c+it)=\Lambda_\eta(c+it)$. Since $f_c''\in L^1(\mathbb R)$, integration by parts twice gives $\widehat f_c(t)=O((1+|t|)^{-2})$, so $\widehat f_c\in L^1(\mathbb R)$. Fourier inversion at $y=\log x$ yields \eqref{eq:theta-Mellin-inversion}.
\end{proof}

\begin{theorem}[Mellin realization and uniqueness of the shifted eta kernel]\label{thm:shifted-kernel-uniqueness}
For any $\alpha\in\mathbb C$ and every $s$ with $\RePart(s+\alpha)>0$,
\begin{equation}\label{eq:shifted-kernel-mellin}
\int_0^\infty x^{s-1}\,x^\alpha K(x)\,dx=\Lambda_\eta(s+\alpha).
\end{equation}
Moreover, fix $c\in\mathbb R$ with $c+\RePart(\alpha)>0$. If measurable kernels $k_1,k_2$ satisfy
\[
\int_0^\infty x^{c-1}|k_j(x)|\,dx<\infty\quad(j=1,2)
\]
and their Mellin transforms agree at every $c+it$, $t\in\mathbb R$, then $k_1=k_2$ almost everywhere. Consequently, if a kernel $k$ in this weighted $L^1$ class has Mellin transform $\Lambda_\eta(s+\alpha)$ on the vertical line $\Re s=c$, then $k(x)=x^\alpha K(x)$ almost everywhere.
\end{theorem}
\begin{proof}
The first assertion follows by replacing $s$ by $s+\alpha$ in Theorem~\ref{thm:theta-eta}; convergence is valid precisely when $\RePart(s+\alpha)>0$. For uniqueness, put $g(y)=e^{cy}(k_1(e^y)-k_2(e^y))$. The substitution $x=e^y$ shows $g\in L^1(\mathbb R)$ and
\[
\widehat g(t):=\int_{\mathbb R}g(y)e^{ity}\,dy
=\int_0^\infty x^{c+it-1}(k_1(x)-k_2(x))\,dx=0.
\]
Uniqueness of the Fourier transform on $L^1(\mathbb R)$ implies $g=0$ almost everywhere. Changing variables back gives $k_1=k_2$ almost everywhere. Apply this result to $k_1=k$ and $k_2(x)=x^\alpha K(x)$; the latter belongs to the required weighted $L^1$ class because $x^{c-1}|x^\alpha K(x)|$ is integrable when $c+\RePart(\alpha)>0$.
\end{proof}

\begin{remark}
For $\alpha=3/8$, the function whose Mellin transform is the shifted completed eta function is exactly $x^{3/8}(1-\theta_4(0\mid ix^2))$. Thus a claim that another kernel has transform $\Lambda_\eta(s+3/8)$ entails, under the stated integrability hypotheses, equality with this weighted Jacobi kernel almost everywhere. This supplies a concrete kernel-level test for the missing correspondence; it does not assert that the source construction passes the test.
\end{remark}

\begin{remark}\label{rem:kernel-not-variant}
Theorem~\ref{thm:theta-eta} is an exact theorem about the standard Jacobi kernel. The weighted series $\ThetaB$ in \eqref{eq:thetaB} is a different object, as shown by Proposition~\ref{prop:thetaB-mellin}. No equality between the manuscript's full variant and $K$ or its Mellin transform is asserted here.
\end{remark}

A second useful representation is the Fermi-type integral
\begin{equation}\label{eq:fermi}
\int_0^\infty\frac{x^{s-1}}{e^x+1}\,dx=\Gamma(s)\etaD(s),\qquad \RePart(s)>0.
\end{equation}
For $\RePart(s)>1$ this follows by expanding $(e^x+1)^{-1}=\sum_{n\ge1}(-1)^{n-1}e^{-nx}$ and integrating term by term. Both sides are holomorphic for $\RePart(s)>0$, so analytic continuation by the identity theorem establishes \eqref{eq:fermi} throughout that half-plane.

\section{The completed eta reflection and the shift \texorpdfstring{$\tfrac14-s$}{one quarter minus s}}\label{sec:shift}

The completed zeta function satisfies the classical meromorphic equation
\begin{equation}\label{eq:completed-zeta-fe}
\LamZ(s)=\LamZ(1-s).
\end{equation}
Since $\LamE(s)=(1-2^{1-s})\LamZ(s)$, we obtain the completed eta equation.
\begin{theorem}[Completed eta reflection]\label{thm:eta-fe}
As a meromorphic identity,
\begin{equation}\label{eq:eta-fe-cross}
\boxed{(1-2^s)\LamE(s)=(1-2^{1-s})\LamE(1-s).}
\end{equation}
\end{theorem}
\begin{proof}
By definition, $\LamE(s)=(1-2^{1-s})\LamZ(s)$ and
\[
\LamE(1-s)=(1-2^s)\LamZ(1-s)=(1-2^s)\LamZ(s).
\]
Multiplication gives \eqref{eq:eta-fe-cross}. The identity is interpreted meromorphically at exceptional points.
\end{proof}

\begin{theorem}[Shifted completed eta equation]\label{thm:shifted-eta}
Let
\begin{equation}\label{eq:F-def}
F(s)=\LamE\!\left(s+\frac38\right).
\end{equation}
Then
\begin{equation}\label{eq:shifted-fe}
\boxed{(1-2^{s+3/8})F(s)=(1-2^{5/8-s})F\!\left(\frac14-s\right)}
\end{equation}
as a meromorphic identity.
\end{theorem}
\begin{proof}
Apply \eqref{eq:eta-fe-cross} at $s+3/8$. The reflected completed eta argument is $1-(s+3/8)=5/8-s$, and
\[
\LamE\!\left(\frac58-s\right)=F\!\left(\frac58-s-\frac38\right)=F\!\left(\frac14-s\right).
\]
Substitution gives \eqref{eq:shifted-fe}.
\end{proof}

\begin{theorem}[The complete family of affine eta reflections]\label{thm:general-shift}
For any $\alpha\in\mathbb C$, define $F_\alpha(s)=\Lambda_\eta(s+\alpha)$ and $c_\alpha=1-2\alpha$. Then
\begin{equation}\label{eq:general-shift-cross}
\boxed{(1-2^{s+\alpha})F_\alpha(s)
=(1-2^{1-s-\alpha})F_\alpha(c_\alpha-s)}
\end{equation}
as a meromorphic identity. Equivalently, wherever division is permitted,
\[
F_\alpha(s)=\frac{1-2^{1-s-\alpha}}{1-2^{s+\alpha}}F_\alpha(c_\alpha-s).
\]
Given any $c\in\mathbb C$, taking $\alpha=(1-c)/2$ produces a reflected argument $c-s$. Moreover, the normalized function
\[
G_\alpha(s)=\frac{F_\alpha(s)}{1-2^{1-s-\alpha}}=\Lambda_\zeta(s+\alpha)
\]
satisfies the pure reflection $G_\alpha(s)=G_\alpha(c_\alpha-s)$ meromorphically.
\end{theorem}
\begin{proof}
Apply Theorem~\ref{thm:eta-fe} with argument $s+\alpha$. The reflected argument of $\Lambda_\eta$ is $1-s-\alpha$, which equals $(c_\alpha-s)+\alpha$. This proves \eqref{eq:general-shift-cross}. The statement about arbitrary $c$ follows by solving $c=1-2\alpha$. Dividing by the elementary factor gives $G_\alpha=\Lambda_\zeta(s+\alpha)$; the pure reflection follows from $\Lambda_\zeta(z)=\Lambda_\zeta(1-z)$.
\end{proof}

\begin{corollary}[The quarter reflection is not diagnostic by itself]\label{cor:quarter-not-diagnostic}
The argument $\tfrac14-s$ is obtained by the choice $\alpha=3/8$, but any other affine reflection $c-s$ is obtained by another shift. Thus observing or postulating the argument $\tfrac14-s$ does not identify a function with the shifted completed eta function unless the shift and the multiplicative normalization are derived independently from the definition of that function.
\end{corollary}

\begin{proposition}[Multipliers preserve the shifted reflection with a correction factor]\label{prop:multiplier-shift}
Let $A$ be meromorphic and set $H_\alpha(s)=A(s)F_\alpha(s)$. Then
\begin{equation}\label{eq:multiplier-shift}
A(c_\alpha-s)(1-2^{s+\alpha})H_\alpha(s)
=A(s)(1-2^{1-s-\alpha})H_\alpha(c_\alpha-s)
\end{equation}
as a meromorphic identity.
\end{proposition}
\begin{proof}
Substitute $F_\alpha(s)=H_\alpha(s)/A(s)$ into \eqref{eq:general-shift-cross} and clear denominators on the open set where they are nonzero. The resulting identity extends meromorphically.
\end{proof}

\begin{remark}[Meaning of the quarter-shift]\label{rem:shift-status}
Theorem~\ref{thm:shifted-eta} establishes that a reflection $s\mapsto\tfrac14-s$ is exact for the shifted completed eta function $F$. It is a reparameterization of the classical eta reflection, not evidence by itself that the original variant has this symmetry. The fixed line of the reflection is $\RePart(s)=1/8$. A functional equation for a logarithm requires, in addition, a proved multiplicative identity for the relevant function and a branch convention on a domain avoiding its zeros and poles.
\end{remark}

\section{A corrected analytic surrogate and its eta factorization}\label{sec:corrected-model}

Because the real-part-based limit in \eqref{eq:Pn} fails for nonreal $s$, an analytic model must choose a valid gamma factor explicitly. The displayed product in the motivating manuscript contains $\Pi(s/2)$; the standard Euler limit \eqref{eq:euler} makes this $\Gamma(1+s/2)$. We use that choice in the following model. This is a corrected interpretation of the displayed product, not a claim that it reproduces every operation of the supplied legacy C routine.

\begin{definition}[Corrected analytic product]\label{def:V}
For $\RePart(s)>0$ away from poles, define
\begin{equation}\label{eq:V}
V(s)=\Gamma\!\left(1+\frac{s}{2}\right)\ThetaB(s)\pi^{-s/2}\zeta(s),\qquad
\ThetaB(s)=\frac{e^{-s}}{(1-e^{-s})^2}.
\end{equation}
The expression defines a meromorphic function by continuation wherever the factors are meromorphically defined.
\end{definition}

\begin{theorem}[Completed-zeta factorization and reflection of the corrected product]\label{thm:V-zeta-FE}
Let $A(s)=\tfrac{s}{2}\ThetaB(s)$ and let $V$ be as in Definition~\ref{def:V}. Then
\begin{equation}\label{eq:V-zeta-factorization}
\boxed{V(s)=A(s)\Lambda_\zeta(s).}
\end{equation}
In particular, $V$ satisfies the meromorphic functional equation
\begin{equation}\label{eq:V-reflection}
\boxed{V(s)=\frac{s\ThetaB(s)}{(1-s)\ThetaB(1-s)}V(1-s)
=\frac{s}{1-s}\frac{\sinh^2((1-s)/2)}{\sinh^2(s/2)}V(1-s).}
\end{equation}
After normalization, $V(s)/A(s)=\Lambda_\zeta(s)$ has the classical pure reflection. Its shifted normalization $V(s+3/8)/A(s+3/8)$ therefore has a pure $s\mapsto\tfrac14-s$ reflection. The unnormalized product $V(s)$ itself is governed by \eqref{eq:V-reflection}, not by a pure quarter reflection.
\end{theorem}
\begin{proof}
By the recurrence $\Gamma(1+s/2)=(s/2)\Gamma(s/2)$, Definition~\ref{def:V} gives
\[
V(s)=\frac{s}{2}\ThetaB(s)\pi^{-s/2}\Gamma(s/2)\zeta(s)=A(s)\Lambda_\zeta(s),
\]
which proves \eqref{eq:V-zeta-factorization}. The completed zeta functional equation gives
\[
V(s)=A(s)\Lambda_\zeta(s)=A(s)\Lambda_\zeta(1-s)=\frac{A(s)}{A(1-s)}V(1-s),
\]
where defined, with meromorphic continuation elsewhere. Substituting $A(s)=s\ThetaB(s)/2$ and \eqref{eq:thetaB-hyperbolic} yields \eqref{eq:V-reflection}. The statements about normalized shifts follow directly by division by $A$.
\end{proof}

\begin{proposition}[Exact eta factorization of the corrected product]\label{prop:V-factorization}
At points where the displayed factors are finite and $1-2^{1-s}\ne0$,
\begin{equation}\label{eq:V-factorization}
\boxed{V(s)=M_0(s)\LamE(s),\qquad M_0(s)=\frac{s\,\ThetaB(s)}{2(1-2^{1-s})}.}
\end{equation}
The identity extends meromorphically.
\end{proposition}
\begin{proof}
Using $\etaD(s)=(1-2^{1-s})\zeta(s)$, we have
\[
\LamE(s)=\pi^{-s/2}\Gamma\!\left(\frac{s}{2}\right)(1-2^{1-s})\zeta(s).
\]
The gamma recurrence $\Gamma(1+s/2)=(s/2)\Gamma(s/2)$ gives \eqref{eq:V-factorization}.
\end{proof}

\begin{remark}[Novelty status]\label{rem:factorization-status}
The factorization \eqref{eq:V-factorization} is exact but is an immediate algebraic consequence of the product definition and the standard relation between eta and zeta. It is not a proof of the Jacobi theta Mellin identity for $V$, and it does not establish a $\tfrac14-s$ reflection for $V$ itself. The normalized shifted function $V(s+3/8)/M_0(s+3/8)$ equals $F(s)$ and therefore obeys Theorem~\ref{thm:shifted-eta}; this is a normalization of the classical eta identity, not a new functional equation for the unnormalized product.
\end{remark}

A key consequence is that at any nontrivial zero $\rho$ of $\zeta$, the corrected product satisfies $V(\rho)=0$, provided its other factors are finite there. Hence $\log V(\rho)$ is undefined. A finite approximation $V_N(\rho)$ may be nonzero, but its logarithm is then a logarithm of a truncation residual and need not converge to a finite value. This is important when interpreting logarithmic plots formed from finite sums at zeta zeros.

\begin{proposition}[Zero preservation and the logarithmic obstruction]\label{prop:zero-log}
Let $\rho$ be a nontrivial zero of $\zeta$. Then the corrected product $V$ is holomorphic at $\rho$ and has a zero of the same multiplicity as $\zeta$. In particular, $\log V(\rho)$ is undefined for every branch of the complex logarithm. A finite truncation $V_N(\rho)$ may be nonzero, but $\log V_N(\rho)$ then measures the logarithm of a truncation residual; it is not a value of $\log V$ at $\rho$.
\end{proposition}
\begin{proof}
A nontrivial zero satisfies $0<\RePart(\rho)<1$, so neither $\Gamma(1+\rho/2)$ nor $\pi^{-\rho/2}$ has a pole there. Also $\ThetaB(\rho)=1/(4\sinh^2(\rho/2))$ is finite and nonzero because its poles lie on $2\pi i\mathbb Z$, which is disjoint from the critical strip. Thus $V(s)=C(s)\zeta(s)$ with $C$ holomorphic and nonzero in a neighborhood of $\rho$. Multiplication by such a factor preserves the zero multiplicity. Since a complex logarithm is defined only at nonzero values, $\log V(\rho)$ does not exist.
\end{proof}

\begin{theorem}[Local singularity structure of the corrected product]\label{thm:V-singularities}
Let $V$ be the meromorphic function in Definition~\ref{def:V}. Its Laurent expansions at $s=0$ and $s=1$ are
\begin{align}
V(s)&=-\frac{1}{2s^2}+\frac{\gamma-\log\pi-2\log2}{4s}+O(1),\qquad s\to0,\label{eq:V-at-zero}\\
V(s)&=\frac{e}{2(e-1)^2}\frac1{s-1}+O(1),\qquad s\to1.\label{eq:V-at-one}
\end{align}
For every integer $n\geq1$, the singularity at $s=-2n$ is removable, and
\begin{equation}\label{eq:V-trivial-values}
V(-2n)=-\frac{(2n)!}{(n-1)!\,2^{2n}\pi^n}\,\Theta^B(-2n)\,\zeta(2n+1)\ne0.
\end{equation}
At a nonzero point $s_0=2\pi i k$, $k\in\mathbb Z\setminus\{0\}$, the order of $V$ is the order of the zero of $\zeta$ at $s_0$ minus two (with order zero for a nonzero value); in particular, it has a double pole there whenever $\zeta(s_0)\ne0$. Apart from the weighted-series pole lattice and $s=1$, the gamma poles are cancelled by the trivial zeros of $\zeta$.
\end{theorem}
\begin{proof}
The expansions
\[
\Theta^B(s)=s^{-2}-\frac1{12}+O(s^2),\quad
\Gamma(1+s/2)=1-\frac{\gamma s}{2}+O(s^2),
\]
\[
\pi^{-s/2}=1-\frac{s\log\pi}{2}+O(s^2),\quad
\zeta(s)=-\frac12-\frac{s\log(2\pi)}2+O(s^2)
\]
follow from the Taylor expansion of $\sinh$, the gamma expansion at $1$, and the classical values $\zeta(0)=-1/2$, $\zeta'(0)=-\tfrac12\log(2\pi)$. Multiplying gives \,\eqref{eq:V-at-zero}. At $s=1$, the zeta function has residue one, while $\Gamma(3/2)\pi^{-1/2}=1/2$ and $\Theta^B(1)=e/(e-1)^2$, proving \,\eqref{eq:V-at-one}.

At $s=-2n$, set $s=-2n+h$. The gamma factor has expansion
\[
\Gamma\!\left(1+\frac{s}{2}\right)=\frac{2(-1)^{n-1}}{(n-1)!}\frac1h+O(1),
\]
while $\zeta(s)=\zeta'(-2n)h+O(h^2)$. The pole and zero cancel. The differentiated zeta functional equation gives
\[
\zeta'(-2n)=\frac{(-1)^n(2n)!}{2(2\pi)^{2n}}\zeta(2n+1),
\]
which yields \,\eqref{eq:V-trivial-values}; nonvanishing follows from $\zeta(2n+1)>0$ and $\Theta^B(-2n)>0$. Finally, at $s_0=2\pi ik\ne0$, the gamma and pi factors are holomorphic and nonzero and $\Theta^B$ has a double pole with nonzero leading coefficient. Multiplication by $\zeta$, which is holomorphic at $s_0$, changes the local order by the order of its zero. The remaining possible gamma poles occur at $-2n$ and were just shown to be removable; the remaining zeta pole is at $1$. This proves the assertion.
\end{proof}

\section{Eta truncation and propagated error bounds}\label{sec:truncation}

For $N\ge1$, set
\begin{equation}\label{eq:etaN}
\etaN(s)=\sum_{n=1}^{N}(-1)^{n-1}n^{-s}.
\end{equation}
\begin{theorem}[Explicit eta tail bound]\label{thm:eta-bound}
If $s=\sigma+it$ with $\sigma>0$, then
\begin{equation}\label{eq:eta-bound}
\boxed{|\etaD(s)-\etaN(s)|\le\frac{|s|}{\sigma}(N+1)^{-\sigma}.}
\end{equation}
Consequently,
\begin{equation}\label{eq:LamE-bound}
\boxed{|\LamE(s)-\pi^{-s/2}\Gamma(s/2)\etaN(s)|
\le \pi^{-\sigma/2}\left|\Gamma\!\left(\frac{s}{2}\right)\right|\frac{|s|}{\sigma}(N+1)^{-\sigma}.}
\end{equation}
The convergence $\etaN\to\etaD$ and the corresponding completed eta convergence are uniform on compact subsets of $\RePart(s)>0$.
\end{theorem}
\begin{proof}
Let $a_n=(-1)^{n-1}$, $b_n=n^{-s}$, and $A_m=\sum_{n=N+1}^m a_n$. Then $|A_m|\le1$. Summation by parts gives, after letting the upper index tend to infinity,
\[
|\etaD(s)-\etaN(s)|\le\sum_{n=N+1}^{\infty}|n^{-s}-(n+1)^{-s}|.
\]
Since
\[
n^{-s}-(n+1)^{-s}=s\int_n^{n+1}x^{-s-1}\,dx,
\]
we obtain
\[
\begin{aligned}
|\etaD(s)-\etaN(s)|
&\le |s|\int_{N+1}^{\infty}x^{-\sigma-1}\,dx\\
&=\frac{|s|}{\sigma}(N+1)^{-\sigma}.
\end{aligned}
\]
Multiplication by $\pi^{-s/2}\Gamma(s/2)$ proves \eqref{eq:LamE-bound}. If $K$ is compact in $\RePart(s)>0$, then $\RePart(s)\ge\delta>0$ and $|s|\le T$ on $K$. The bound is uniform on $K$; the gamma and pi factors are bounded there, proving the local uniform convergence assertions.
\end{proof}

\begin{corollary}[Eta-based zeta approximation]\label{cor:zeta-bound}
At points where $D(s)=1-2^{1-s}\ne0$, define
\begin{equation}\label{eq:zetaNeta}
\widehat\zeta_N(s)=\frac{\etaN(s)}{1-2^{1-s}}.
\end{equation}
Then, for $\RePart(s)=\sigma>0$,
\begin{equation}\label{eq:zeta-error-bound}
\boxed{|\zeta(s)-\widehat\zeta_N(s)|\le
\frac{|s|}{\sigma|1-2^{1-s}|}(N+1)^{-\sigma}.}
\end{equation}
The approximation converges locally uniformly on compact subsets of $\{s:\RePart(s)>0,\ D(s)\ne0\}$.
\end{corollary}
\begin{proof}
Divide \eqref{eq:eta-bound} by $|D(s)|$. On a compact set avoiding the zero set of $D$, the denominator is bounded away from zero, giving local uniform convergence.
\end{proof}

\begin{corollary}[Uniform zeta error away from $\RePart(s)=1$]\label{cor:zeta-uniform-strip}
Fix $0<\delta<1/2$ and $T>0$. For $s=\sigma+it$ satisfying $\delta\leq\sigma\leq1-\delta$ and $|t|\leq T$, one has
\begin{equation}\label{eq:zeta-uniform-strip}
\left|\zeta(s)-\widehat\zeta_N(s)\right|
\leq \frac{\sqrt{(1-\delta)^2+T^2}}{\delta(2^\delta-1)}(N+1)^{-\delta}.
\end{equation}
In particular, on the critical line $\sigma=1/2$,
\begin{equation}\label{eq:zeta-critical-line-bound}
\left|\zeta\!\left(\frac12+it\right)-\widehat\zeta_N\!\left(\frac12+it\right)\right|
\leq \frac{2|\frac12+it|}{\sqrt2-1}(N+1)^{-1/2}.
\end{equation}
\end{corollary}
\begin{proof}
The reverse triangle inequality gives
\[
|1-2^{1-s}|\geq\bigl|2^{1-\sigma}-1\bigr|.
\]
When $\delta\leq\sigma\leq1-\delta$, this is at least $2^\delta-1$, while $|s|\leq\sqrt{(1-\delta)^2+T^2}$ and $\sigma\geq\delta$. Substitution in Corollary~\ref{cor:zeta-bound}, together with $(N+1)^{-\sigma}\leq(N+1)^{-\delta}$, proves \eqref{eq:zeta-uniform-strip}. If $\sigma=1/2$, then $|2^{1-s}|=\sqrt2$, so $|1-2^{1-s}|\geq\sqrt2-1$; the general bound yields \eqref{eq:zeta-critical-line-bound}.
\end{proof}

\begin{corollary}[Error bound for the corrected product]\label{cor:V-bound}
For $V$ in Definition~\ref{def:V}, define
\[
V_N(s)=\Gamma\!\left(1+\frac{s}{2}\right)\ThetaB(s)\pi^{-s/2}\widehat\zeta_N(s).
\]
Then
\begin{equation}\label{eq:V-bound}
\boxed{|V(s)-V_N(s)|\le
\pi^{-\sigma/2}\left|\Gamma\!\left(1+\frac{s}{2}\right)\right|\,|\ThetaB(s)|
\frac{|s|}{\sigma|1-2^{1-s}|}(N+1)^{-\sigma}.}
\end{equation}
\end{corollary}
\begin{proof}
The factor multiplying $\zeta(s)$ in \eqref{eq:V} is independent of $N$. Multiplying the estimate in Corollary~\ref{cor:zeta-bound} by its modulus proves the result.
\end{proof}

\begin{theorem}[Error bound when both product factors are truncated]\label{thm:both-truncated}
Fix $s=\sigma+it$ with $\sigma>0$ and $D(s)=1-2^{1-s}\ne0$. Let
\[
\widehat\zeta_N(s)=\frac{\eta_N(s)}{D(s)},\qquad
V_N(s)=\Gamma\!\left(1+\frac{s}{2}\right)\pi^{-s/2}\Theta^B_N(s)\widehat\zeta_N(s),
\]
where $\Theta^B_N$ is the finite sum in Proposition~\ref{prop:thetaB-tail}. Put $r=e^{-\sigma}$ and
\[
R_N(\sigma)=r^{N+1}\left(\frac{N+1}{1-r}+\frac{r}{(1-r)^2}\right).
\]
Then
\begin{align}
|V(s)-V_N(s)|\leq &\;\pi^{-\sigma/2}\left|\Gamma\!\left(1+\frac{s}{2}\right)\right|\Bigg[
|\Theta^B_N(s)|\frac{|s|}{\sigma|D(s)|}(N+1)^{-\sigma}
+|\zeta(s)|R_N(\sigma)\Bigg].\label{eq:both-truncated-bound}
\end{align}
Thus $V_N(s)\to V(s)$ for each such $s$. The convergence is locally uniform on compact subsets of $\{s:\RePart(s)>0,\ D(s)\ne0\}$.
\end{theorem}
\begin{proof}
Write $C(s)=\Gamma(1+s/2)\pi^{-s/2}$. The exact decomposition
\[
\ThetaB(s)\zeta(s)-\Theta^B_N(s)\widehat\zeta_N(s)
=\Theta^B_N(s)(\zeta(s)-\widehat\zeta_N(s))
+(\ThetaB(s)-\Theta^B_N(s))\zeta(s)
\]
and the triangle inequality give
\[
|V-V_N|\leq |C(s)|\left(|\Theta^B_N(s)|\,|\zeta-\widehat\zeta_N|+|\zeta(s)|\,|\ThetaB-\Theta^B_N|\right).
\]
Apply Corollary~\ref{cor:zeta-bound} to the first error and Proposition~\ref{prop:thetaB-tail} to the second. This proves \eqref{eq:both-truncated-bound}. On compact sets in the stated domain, $\RePart(s)$ is bounded below by a positive constant and $|D(s)|$ is bounded away from zero; the factors are bounded and both tail bounds tend uniformly to zero.
\end{proof}

\subsection{A geometrically convergent Euler transform for eta}
The ordinary eta partial sums have algebraic convergence on the critical line. A second representation follows by applying the Euler transformation for alternating series to the eta series. For $k\ge0$, set
\begin{equation}\label{eq:euler-difference}
\Delta_k(s)=\sum_{j=0}^{k}(-1)^j\binom{k}{j}(j+1)^{-s},\qquad
\eta^{\mathrm E}_K(s)=\sum_{k=0}^{K}\frac{\Delta_k(s)}{2^{k+1}}.
\end{equation}
\begin{theorem}[Euler-transformed eta series and explicit geometric error]\label{thm:euler-eta}
For every $s\in\mathbb C$ with $\sigma=\RePart(s)>0$,
\begin{equation}\label{eq:euler-eta}
\boxed{\eta(s)=\sum_{k=0}^{\infty}\frac{\Delta_k(s)}{2^{k+1}}.}
\end{equation}
The series converges absolutely, and its truncation satisfies
\begin{equation}\label{eq:euler-eta-bound}
\boxed{\left|\eta(s)-\eta^{\mathrm E}_K(s)\right|
\leq\frac{\Gamma(\sigma)}{|\Gamma(s)|}\,2^{-(K+1)}.}
\end{equation}
The convergence is locally uniform on $\RePart(s)>0$. For every compact $Q$ in this half-plane, the bound is uniform on $Q$ with constant $\sup_{s\in Q}\Gamma(\RePart(s))/|\Gamma(s)|<\infty$.
\end{theorem}
\begin{proof}
The Laplace representation $ (j+1)^{-s}=\Gamma(s)^{-1}\int_0^\infty u^{s-1}e^{-(j+1)u}\,du$ gives
\begin{equation}\label{eq:euler-difference-integral}
\Delta_k(s)=\frac1{\Gamma(s)}\int_0^\infty u^{s-1}e^{-u}(1-e^{-u})^k\,du.
\end{equation}
Because $0\leq1-e^{-u}<1$, Tonelli's theorem applied to absolute values shows
\[
\sum_{k=0}^{\infty}\frac{|\Delta_k(s)|}{2^{k+1}}
\leq\frac1{|\Gamma(s)|}\int_0^\infty u^{\sigma-1}e^{-u}
\sum_{k=0}^{\infty}\frac{(1-e^{-u})^k}{2^{k+1}}\,du
=\frac1{|\Gamma(s)|}\int_0^\infty\frac{u^{\sigma-1}}{e^u+1}\,du<\infty.
\]
Thus the transformed series is absolutely convergent and its terms may be interchanged with the integral. Summing the geometric series gives
\[
\sum_{k=0}^{\infty}\frac{\Delta_k(s)}{2^{k+1}}
=\frac1{\Gamma(s)}\int_0^\infty\frac{u^{s-1}}{e^u+1}\,du=\eta(s),
\]
where the last identity is the Fermi-type representation \,\eqref{eq:fermi}.

For the tail, summing the geometric series from $k=K+1$ yields
\[
\sum_{k=K+1}^{\infty}\frac{(1-e^{-u})^k}{2^{k+1}}
=\frac{(1-e^{-u})^{K+1}}{2^{K+1}(1+e^{-u})}.
\]
Consequently,
\[
\begin{aligned}
|\eta(s)-\eta_K^{\mathrm E}(s)|
&\leq\frac1{|\Gamma(s)|2^{K+1}}
\int_0^\infty u^{\sigma-1}\frac{e^{-u}(1-e^{-u})^{K+1}}{1+e^{-u}}\,du\\
&\leq\frac1{|\Gamma(s)|2^{K+1}}\int_0^\infty u^{\sigma-1}e^{-u}\,du
=\frac{\Gamma(\sigma)}{|\Gamma(s)|}2^{-(K+1)}.
\end{aligned}
\]
The gamma function has no zeros or poles in $\RePart(s)>0$, so the constant is continuous and bounded on any compact subset; the same estimate proves local uniform convergence.
\end{proof}
\begin{remark}[Numerical stability]
The outer series has geometric truncation error, but the finite differences in \,\eqref{eq:euler-difference} involve alternating binomial sums and can suffer severe cancellation in ordinary precision. A practical implementation should use sufficiently high precision or a stable integral/recurrence evaluation. The bound in \,\eqref{eq:euler-eta-bound} is rigorous but may be conservative for large $|\Im(s)|$ because $|\Gamma(s)|$ is small.
\end{remark}

\subsection{A sharper alternating-tail asymptotic}
The upper bound in Theorem~\ref{thm:eta-bound} is uniform and convenient, but it is not asymptotically sharp. A two-term expansion, with a rigorous remainder, explains the observed slope.
\begin{theorem}[Alternating eta-tail expansion]\label{thm:eta-tail-expansion}
Fix $s\in\mathbb C$ with $\sigma=\RePart(s)>0$ and put $m=N+1$. As $N\to\infty$,
\begin{equation}\label{eq:eta-tail-expansion}
\eta(s)-\eta_N(s)=(-1)^N\left[\frac12m^{-s}+\frac{s}{4}m^{-s-1}-\frac{s(s+1)(s+2)}{48}m^{-s-3}+O_s(m^{-\sigma-5})\right].
\end{equation}
In particular,
\begin{equation}\label{eq:eta-tail-two-term}
\eta(s)-\eta_N(s)=\frac{(-1)^N}{2}m^{-s}\left(1+\frac{s}{2m}+O_s(m^{-3})\right).
\end{equation}
If $\rho$ is a zero of $\zeta$ with $0<\RePart(\rho)=\beta<1$, then $\eta(\rho)=0$ and
\begin{equation}\label{eq:eta-at-zero-asymptotic}
|\eta_N(\rho)|=\frac12m^{-\beta}\left(1+\frac{\beta}{2m}+O_\rho(m^{-2})\right).
\end{equation}
Consequently,
\begin{equation}\label{eq:log-residual-asymptotic}
\log|\eta_N(\rho)|=-\beta\log m-\log2+\frac{\beta}{2m}+O_\rho(m^{-2}),
\end{equation}
and the asymptotic log--log slope is $-\beta$. For fixed $q>1$, the secant slope between indices $N$ and $\lfloor qN\rfloor$ is $-\beta+O_{\rho,q}(N^{-1})$.
\end{theorem}
\begin{proof}
Let $m=N+1$. For $\RePart(s)>0$, the alternating tail admits the Laplace representation
\begin{equation}\label{eq:tail-laplace}
\sum_{k=0}^{\infty}(-1)^k(m+k)^{-s}
=\frac1{\Gamma(s)}\int_0^\infty\frac{u^{s-1}e^{-mu}}{1+e^{-u}}\,du.
\end{equation}
To justify it, first integrate the finite sum through $k=M$ using
\[
(m+k)^{-s}=\frac1{\Gamma(s)}\int_0^\infty u^{s-1}e^{-(m+k)u}\,du.
\]
The finite geometric sum is $[1-(-e^{-u})^{M+1}]/(1+e^{-u})$. The remaining integral is bounded in modulus by
\[
\frac1{|\Gamma(s)|}\int_0^\infty u^{\sigma-1}e^{-(m+M+1)u}\,du
=\frac{\Gamma(\sigma)}{|\Gamma(s)|}(m+M+1)^{-\sigma},
\]
which tends to zero. The limit gives \eqref{eq:tail-laplace}.

Near $u=0$,
\[
\frac1{1+e^{-u}}=\frac12+\frac{u}{4}-\frac{u^3}{48}+O(u^5).
\]
Insert this expansion into \eqref{eq:tail-laplace}. Integration of each monomial gives
$\Gamma(s+j)m^{-s-j}/\Gamma(s)$, and hence the three displayed terms in \eqref{eq:eta-tail-expansion}. To bound the remainder, split the integral at $u=1$. On $0\leq u\leq1$, the Taylor remainder is bounded by $C u^5$; its integral is $O_s(m^{-\sigma-5})$. On $u\geq1$, the difference between the bounded function $(1+e^{-u})^{-1}$ and the degree-three polynomial has at most polynomial growth, while $e^{-mu}$ makes its integral exponentially small in $m$. This proves \eqref{eq:eta-tail-expansion}, and factoring out $\tfrac12m^{-s}$ yields \eqref{eq:eta-tail-two-term}.

At a zeta zero $\rho$ in the critical strip, $1-2^{1-\rho}\ne0$, so $\eta(\rho)=0$. Thus $|\eta_N(\rho)|$ equals the modulus of the tail. Taking the modulus of \eqref{eq:eta-tail-two-term}, and using
$|1+\rho/(2m)+O_\rho(m^{-3})|=1+\beta/(2m)+O_\rho(m^{-2})$,
gives \eqref{eq:eta-at-zero-asymptotic}. Taking logarithms yields \eqref{eq:log-residual-asymptotic}. Applying that expansion at $N$ and $\lfloor qN\rfloor$ proves the final secant-slope statement.
\end{proof}

\begin{theorem}[Finite-product residual at a zeta zero]\label{thm:product-zero-residual}
Let $\rho=\beta+i\gamma$ be a nontrivial zero of $\zeta$, with $0<\beta<1$, and let $V_N$ be the doubly truncated product in Theorem~\ref{thm:both-truncated}. Define
\[
D(\rho)=1-2^{1-\rho},\qquad
K_\rho=\frac{\Gamma(1+\rho/2)\pi^{-\rho/2}\ThetaB(\rho)}{2D(\rho)}.
\]
Then $K_\rho\ne0$ and, as $N\to\infty$,
\begin{equation}\label{eq:product-zero-asymptotic}
V_N(\rho)=(-1)^{N+1}K_\rho(N+1)^{-\rho}
\left(1+\frac{\rho}{2(N+1)}+O_\rho(N^{-2})\right).
\end{equation}
Consequently,
\begin{align}
\log|V_N(\rho)|&=-\beta\log(N+1)+\log|K_\rho|+\frac{\beta}{2(N+1)}+O_\rho(N^{-2}),\label{eq:product-zero-log}\\
\end{align}
and the log--log slope of the truncation residual tends to $-\beta$. In particular, at a zero on the critical line the asymptotic slope is $-1/2$. Along either fixed-parity subsequence of sufficiently large $N$, a consistently continued logarithm also satisfies
\[
\log V_N(\rho)=C_{\rho,\varepsilon}-\rho\log(N+1)+\frac{\rho}{2(N+1)}+O_\rho(N^{-2}),
\]
where $\varepsilon=N\bmod2$ and $C_{\rho,\varepsilon}$ is constant on that subsequence. The alternating sign prevents one from assigning a single unwrapped phase formula across both parities without an additional convention.
\end{theorem}
\begin{proof}
Since $0<\beta<1$, $D(\rho)\ne0$. The factor $\Gamma(1+\rho/2)$ is finite and nonzero, and $\ThetaB(\rho)$ is finite and nonzero because its poles lie on the imaginary axis. Hence $K_\rho\ne0$. By Theorem~\ref{thm:eta-tail-expansion} and $\eta(\rho)=0$,
\[
\eta_N(\rho)=(-1)^{N+1}\frac12(N+1)^{-\rho}\left(1+\frac{\rho}{2(N+1)}+O_\rho(N^{-3})\right).
\]
Also Proposition~\ref{prop:thetaB-tail} gives $\Theta^B_N(\rho)=\ThetaB(\rho)(1+O_\rho(Ne^{-\beta N}))$, whose relative error is $O_\rho(N^{-2})$ (indeed it is smaller than every inverse power). Substituting these formulas into
\[
V_N(\rho)=\Gamma\!\left(1+\frac{\rho}{2}\right)\pi^{-\rho/2}\Theta^B_N(\rho)\frac{\eta_N(\rho)}{D(\rho)}
\]
proves \eqref{eq:product-zero-asymptotic}. Taking moduli and using
$|1+\rho/(2m)+O(m^{-2})|=1+\beta/(2m)+O(m^{-2})$ proves \eqref{eq:product-zero-log}. Taking a logarithm of the factor close to $1$ gives the parity-subsequence formula. Since $V(\rho)=0$, the quantity $V_N(\rho)$ is a truncation residual, not a value of $V$ at a point where $\log V$ exists.
\end{proof}

\section{Numerical experiments and reproducibility}\label{sec:numerics}

\subsection{What is recomputed}

The original plots are not reproduced as if they were validated outputs: the supplied source refers to an external zero list (\texttt{zero1.h}) and original figure files/raw samples that were not supplied as standalone inputs, and the legacy definitions and logarithm implementation require correction. Instead, this section reruns a fully specified corrected calculation based on Definition~\ref{def:V} and Corollaries~\ref{cor:zeta-bound}--\ref{cor:V-bound}. The code uses eta partial sums to evaluate $\zeta$ in the critical strip, truncates the weighted exponential series using its exact geometric tail, and records observed errors together with rigorous truncation bounds. This is a reproducible replacement experiment, not a claim to have recreated the old crossing-index regressions.

Reference special-function values and ordinary sums were computed at 50 decimal digits with \texttt{mpmath}; the Euler-transformed binomial sums were computed at 100 decimal digits to control cancellation. Compensated summation was used for the long ordinary finite sums. The truncation bounds below are analytic bounds for exact sums. Floating-point or arbitrary-precision roundoff is a separate source of numerical error and is not included in those analytic tail bounds. The attached script \texttt{reproduce\_experiments.py} regenerates all tables, CSV files, and figures.

\subsection{Eta truncation at the first nontrivial zero}

Let $\rho_1=\tfrac12+14.134725141734693790\ldots i$ be the first nontrivial zeta zero. The table reports the absolute eta-series truncation error and the bound \eqref{eq:eta-bound} for selected powers of two.

\begin{table}[ht]
\centering
\small
\caption{Eta truncation at $\rho_1$. The error is $|\eta(\rho_1)-\eta_N(\rho_1)|$.}
\label{tab:eta-errors}
\begin{tabular}{r r r r}
\toprule
$N$ & absolute error & bound & error/bound\\
\midrule
$256$ & $3.12313\times10^{-2}$ & $1.76450$ & $1.76998\times10^{-2}$\\
$1{,}024$ & $1.56216\times10^{-2}$ & $8.83542\times10^{-1}$ & $1.76806\times10^{-2}$\\
$4{,}096$ & $7.81203\times10^{-3}$ & $4.41932\times10^{-1}$ & $1.76770\times10^{-2}$\\
$16{,}384$ & $3.90619\times10^{-3}$ & $2.20986\times10^{-1}$ & $1.76762\times10^{-2}$\\
$65{,}536$ & $1.95312\times10^{-3}$ & $1.10496\times10^{-1}$ & $1.76759\times10^{-2}$\\
$262{,}144$ & $9.76562\times10^{-4}$ & $5.52482\times10^{-2}$ & $1.76759\times10^{-2}$\\
\bottomrule
\end{tabular}
\end{table}

\begin{figure}[ht]
\centering
\includegraphics[width=0.86\linewidth]{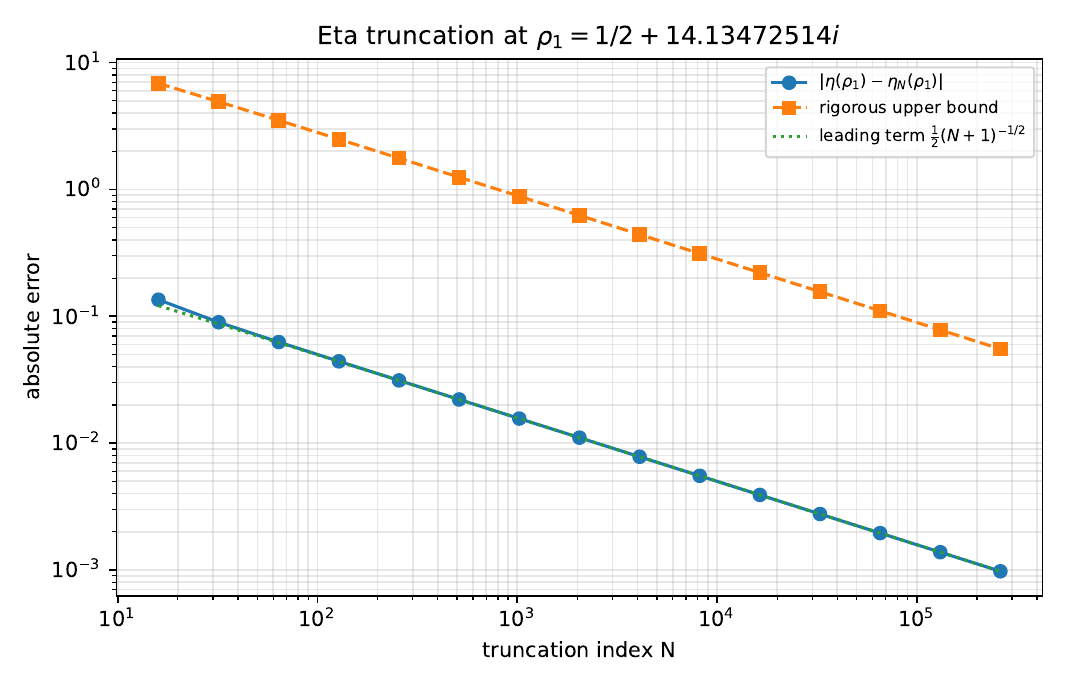}
\caption{Absolute eta truncation error and the explicit upper bound at the first nontrivial zero. Both decay at the scale $N^{-1/2}$; the general bound is conservative.}
\label{fig:eta-error}
\end{figure}

A least-squares fit of $\log_{10}|\eta(\rho_1)-\eta_N(\rho_1)|$ against $\log_{10}N$, for 11 powers of two from $N=256$ through $N=262144$, gives slope $-0.4999284$ and $R^2=0.9999999911$. This is a deterministic diagnostic fit, not an inferential regression; no statistical confidence interval is attached. Its value agrees with the exponent predicted by \eqref{eq:eta-at-zero-asymptotic}.

\begin{figure}[ht]
\centering
\includegraphics[width=0.86\linewidth]{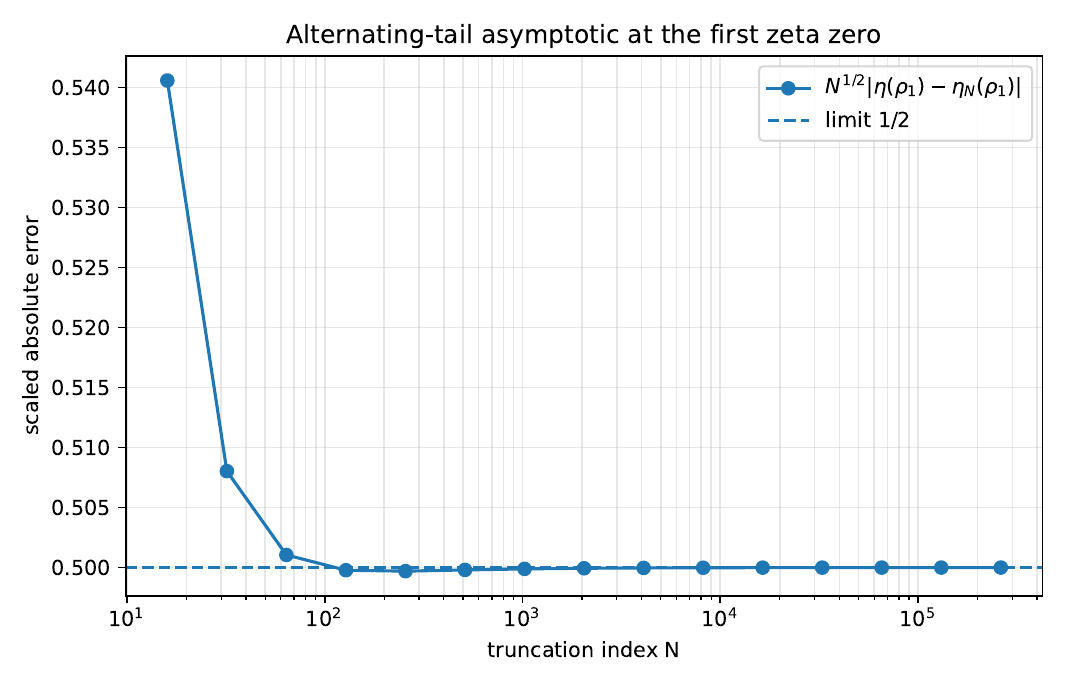}
\caption{The scaled residual $N^{1/2}|\eta(\rho_1)-\eta_N(\rho_1)|$ tends toward $1/2$, as predicted by the first-term alternating-tail asymptotic.}
\label{fig:eta-scaled}
\end{figure}

The same asymptotic propagates to the corrected product. The following table compares $|V_N(\rho_1)|$ with its leading asymptotic. The complete data, including the rigorous bound, are in \texttt{corrected\_product\_zero\_residual.csv}. For $N=262144$, the ratio to the leading asymptotic is approximately $1.0000009540$; the log--log regression over $N=256,512,\ldots,262144$ has slope $-0.4999284022$ and $R^2=0.9999999911$. The slope is predicted by Theorem~\ref{thm:product-zero-residual} and is not evidence for a functional equation.
\begin{table}[ht]
\centering
\small
\caption{Corrected-product residual at the first nontrivial zero and its leading asymptotic.}
\label{tab:product-zero-residual}
\begin{tabular}{r r r r}
\toprule
$N$ & $|V_N(\rho_1)|$ & leading term & ratio \\
\midrule
$256$ & $7.20662\times10^{-7}$ & $7.19688\times10^{-7}$ & $1.00135282$ \\
$1{,}024$ & $3.60467\times10^{-7}$ & $3.60370\times10^{-7}$ & $1.00026770$ \\
$4{,}096$ & $1.80262\times10^{-7}$ & $1.80251\times10^{-7}$ & $1.00006251$ \\
$16{,}384$ & $9.01352\times10^{-8}$ & $9.01338\times10^{-8}$ & $1.00001535$ \\
$65{,}536$ & $4.50681\times10^{-8}$ & $4.50679\times10^{-8}$ & $1.00000382$ \\
$262{,}144$ & $2.25341\times10^{-8}$ & $2.25341\times10^{-8}$ & $1.00000095$ \\
\bottomrule
\end{tabular}
\end{table}
\begin{figure}[ht]
\centering
\includegraphics[width=0.86\linewidth]{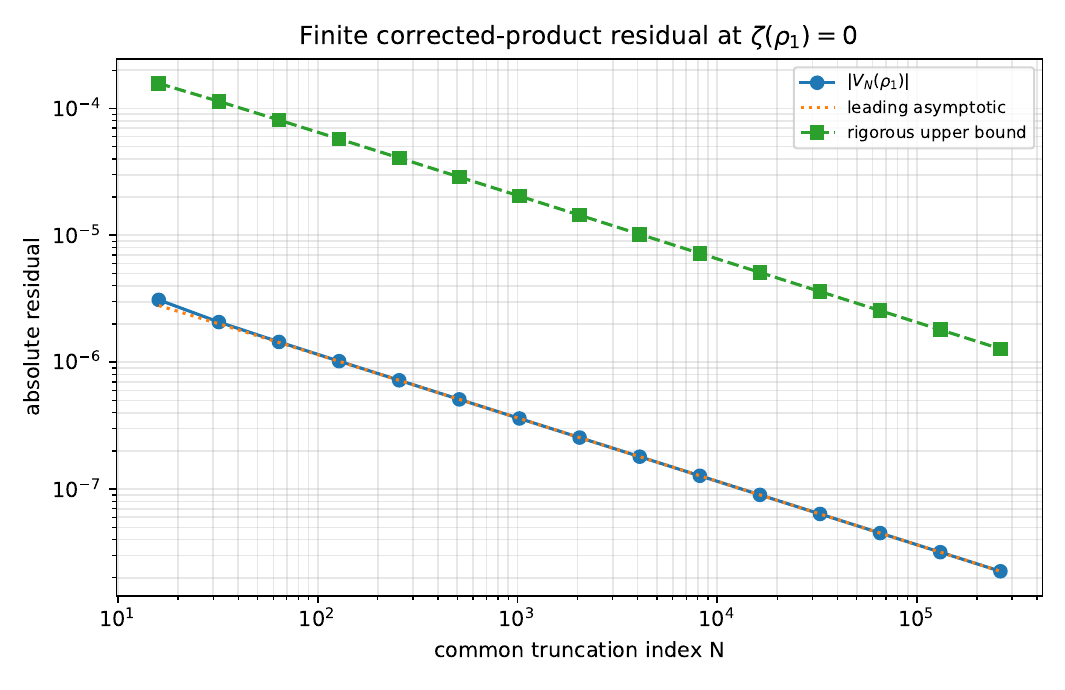}
\caption{Finite corrected-product residual at the first nontrivial zeta zero. The residual follows the leading $(N+1)^{-1/2}$ asymptotic; the rigorous bound is conservative.}
\label{fig:product-zero-residual}
\end{figure}

\subsection{Euler-transformed eta acceleration}
We also evaluate the Euler-transformed series at the first nontrivial zero. The finite-difference sums are computed at 100 decimal digits to mitigate binomial cancellation. The limiting eta value is zero at an exact zeta zero; the ordinate supplied to the calculation is a 100-digit numerical approximation, whose root-location error is far below the displayed residuals. Selected results are shown in Table~\ref{tab:euler-eta} and the full data are included in \texttt{eta\_euler\_acceleration.csv}. For example, the residual at $K=48$ is approximately $9.08\times10^{-11}$, substantially smaller than the $\approx9.77\times10^{-4}$ residual of the ordinary partial sum with $N=262144$. The number of scalar summands in the finite-difference definition through $K=48$ is $\sum_{k=0}^{48}(k+1)=1225$, though its higher arithmetic precision has to be accounted for in a fair cost comparison.

\begin{table}[ht]
\centering
\small
\caption{Euler-transformed eta approximation at the first nontrivial zero; the bound is Theorem~\ref{thm:euler-eta}.}
\label{tab:euler-eta}
\begin{tabular}{r r r r}
\toprule
$K$ & absolute residual & rigorous bound & residual/bound \\
\midrule
$8$ & $5.69937\times10^{-1}$ & $6.06398\times10^{6}$ & $9.39874\times10^{-8}$ \\
$16$ & $2.58476\times10^{-2}$ & $2.36874\times10^{4}$ & $1.09120\times10^{-6}$ \\
$24$ & $3.51790\times10^{-4}$ & $9.25289\times10^{1}$ & $3.80195\times10^{-6}$ \\
$32$ & $2.80445\times10^{-6}$ & $3.61441\times10^{-1}$ & $7.75909\times10^{-6}$ \\
$40$ & $1.71678\times10^{-8}$ & $1.41188\times10^{-3}$ & $1.21595\times10^{-5}$ \\
$48$ & $9.08191\times10^{-11}$ & $5.51515\times10^{-6}$ & $1.64672\times10^{-5}$ \\
$56$ & $4.39925\times10^{-13}$ & $2.15436\times10^{-8}$ & $2.04203\times10^{-5}$ \\
$60$ & $2.99322\times10^{-14}$ & $1.34647\times10^{-9}$ & $2.22301\times10^{-5}$ \\
\bottomrule
\end{tabular}
\end{table}
\begin{figure}[ht]
\centering
\includegraphics[width=0.86\linewidth]{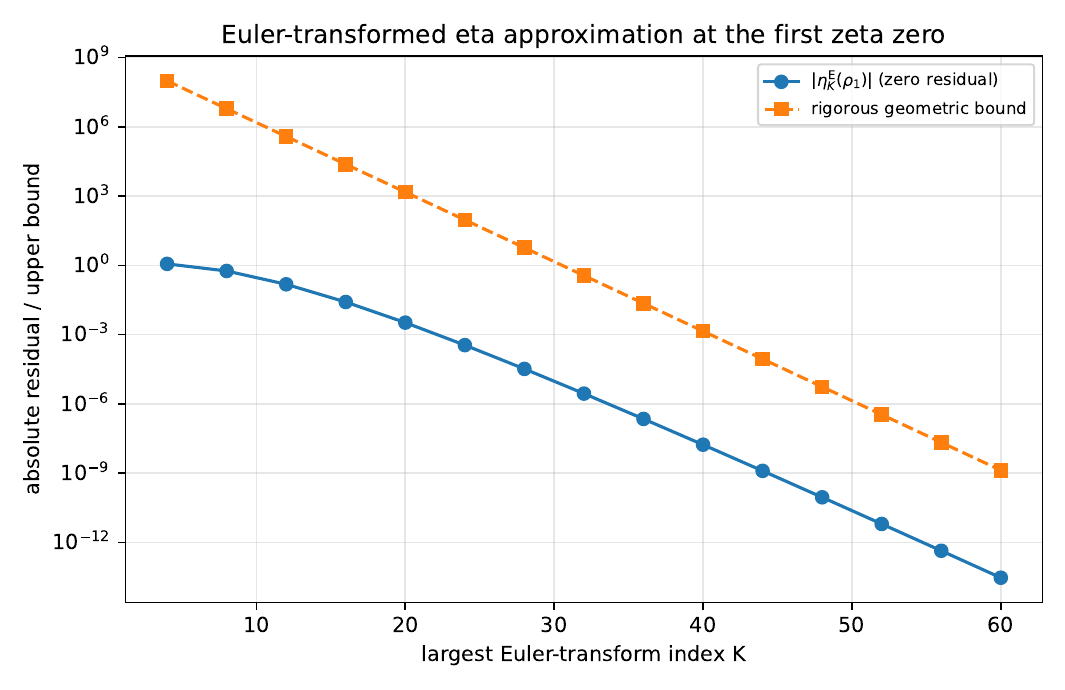}
\caption{Euler-transformed eta residual at the first nontrivial zero and the explicit geometric error bound. The bound is rigorous but conservative for this ordinate.}
\label{fig:euler-eta}
\end{figure}
This experiment demonstrates an alternative convergent representation, not a new functional equation. At large imaginary parts the bound's factor $\Gamma(\sigma)/|\Gamma(s)|$ can be large, so rigorous certification and numerical efficiency must be assessed together.

\subsection{Eta-based evaluation of the zeta function away from a zero}

At $s_1=1/2+82i$, use $\widehat\zeta_N(s_1)=\eta_N(s_1)/(1-2^{1-s_1})$. The propagated bound is given by \eqref{eq:zeta-error-bound}.

\begin{table}[ht]
\centering
\small
\caption{Absolute error and bound for the eta-based zeta approximation at $s_1=1/2+82i$.}
\label{tab:zeta-errors}
\begin{tabular}{r r r r}
\toprule
$N$ & absolute error & bound & error/bound\\
\midrule
$256$ & $5.88035\times10^{-2}$ & $1.90235\times10^{1}$ & $3.09109\times10^{-3}$\\
$1{,}024$ & $2.90715\times10^{-2}$ & $9.52568$ & $3.05191\times10^{-3}$\\
$4{,}096$ & $1.45275\times10^{-2}$ & $4.76458$ & $3.04906\times10^{-3}$\\
$16{,}384$ & $7.26375\times10^{-3}$ & $2.38251$ & $3.04878\times10^{-3}$\\
$65{,}536$ & $3.63191\times10^{-3}$ & $1.19128$ & $3.04874\times10^{-3}$\\
$262{,}144$ & $1.81596\times10^{-3}$ & $5.95645\times10^{-1}$ & $3.04873\times10^{-3}$\\
\bottomrule
\end{tabular}
\end{table}

\begin{figure}[ht]
\centering
\includegraphics[width=0.86\linewidth]{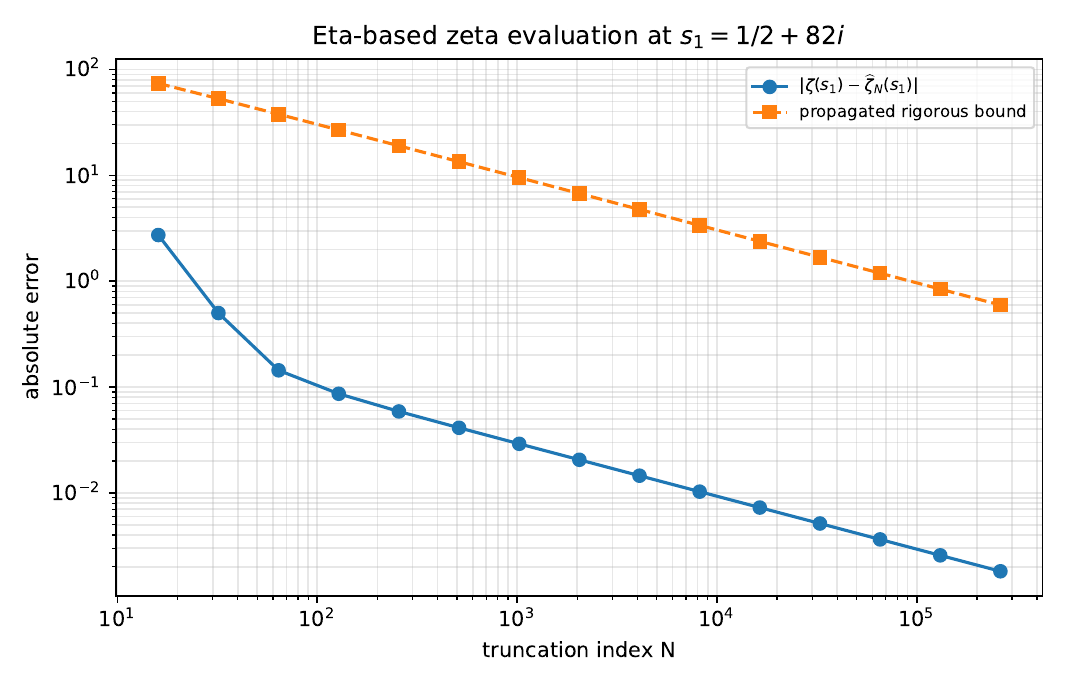}
\caption{Eta-based zeta approximation error and the propagated rigorous truncation bound at $s_1=1/2+82i$.}
\label{fig:zeta-error}
\end{figure}

For 11 powers of two from $N=256$ through $N=262144$, the descriptive least-squares slope of the log absolute zeta error versus $\log_{10}N$ is $-0.5009441$, with $R^2=0.9999939531$. The exact value used for comparison is $\zeta(1/2+82i)=1.229042022307400439-2.4042436204853256541i$ to the displayed precision. The analytic product prefactor
\[
C(s)=\Gamma\!\left(1+\frac{s}{2}\right)\ThetaB(s)\pi^{-s/2}
\]
has $|C(s_1)|\approx9.19097188266\times10^{-27}$. At $N=262144$, the computed error in the fully truncated $V_N(s_1)$ is approximately $1.66904\times10^{-29}$, while the combined bound from \eqref{eq:both-truncated-bound} is approximately $5.47455\times10^{-27}$.

\subsection{Truncating the weighted exponential factor and the full product}
The exact formula and bound in Proposition~\ref{prop:thetaB-tail} were evaluated at $s_B=1/2+14.13472514173469i$. The observed truncation error is below the explicit bound at every tested power of two. The series tail is exponentially small in $N$, in contrast with the algebraic $N^{-1/2}$ eta-series error on the critical line.
\begin{figure}[ht]
\centering
\includegraphics[width=0.86\linewidth]{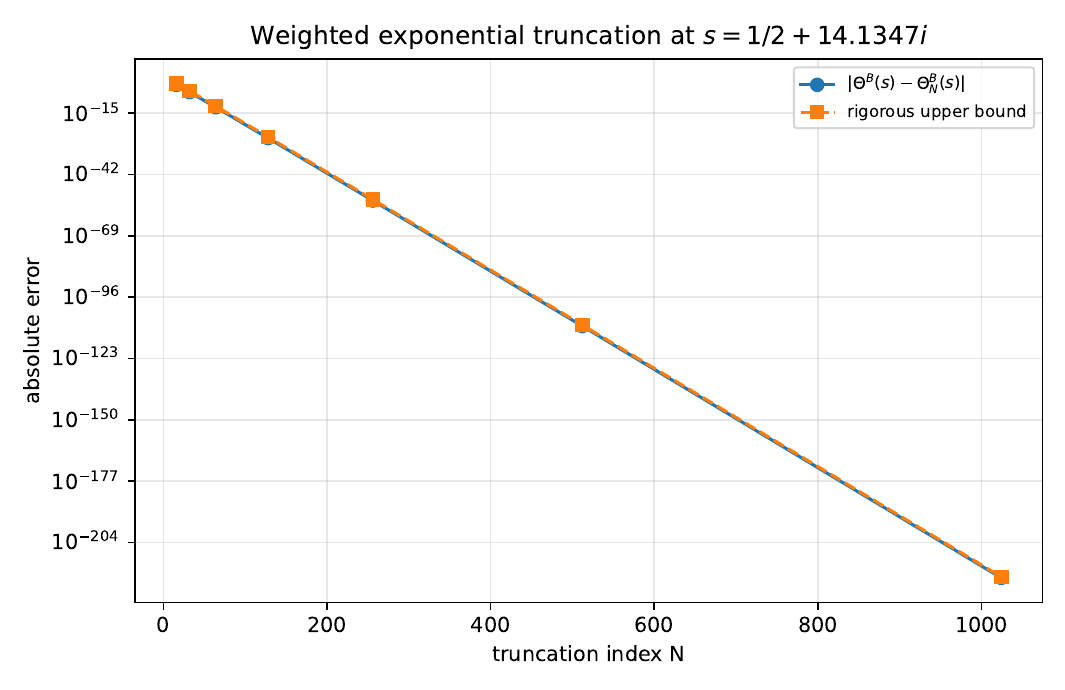}
\caption{Weighted-exponential truncation error and its explicit geometric-tail bound. The exact tail formula was used for evaluation to avoid artificial cancellation from long direct summation.}
\label{fig:weighted-tail}
\end{figure}

For the full finite surrogate, both approximations are truncated at the same index:
\[
V_N(s)=\Gamma\!\left(1+\frac{s}{2}\right)\pi^{-s/2}\Theta^B_N(s)\widehat\zeta_N(s).
\]
The combined error bound is Theorem~\ref{thm:both-truncated}. At $s_1=1/2+82i$ and $N=262144$, the computed error is approximately $1.66904\times10^{-29}$ and the combined rigorous bound is approximately $5.47455\times10^{-27}$. The bound is conservative, as expected; it is an analytical guarantee, not a confidence interval.
\begin{figure}[ht]
\centering
\includegraphics[width=0.86\linewidth]{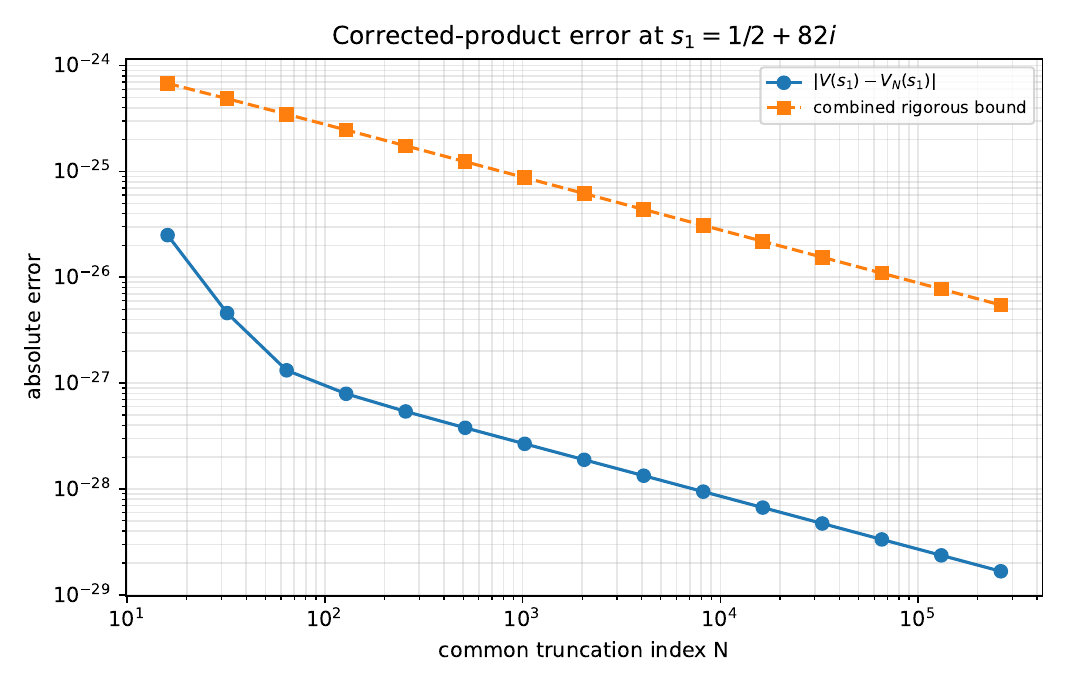}
\caption{Error in the corrected finite product when both the eta-based zeta approximation and the weighted-exponential factor are truncated, with the combined rigorous bound from Theorem~\ref{thm:both-truncated}.}
\label{fig:full-product-error}
\end{figure}

\subsection{High-precision identity checks}

At $s=0.31+7.3i$, 50-decimal-digit evaluation gives residuals of approximately $5.22\times10^{-53}$ for the eta factorization \eqref{eq:V-factorization}, $8.35\times10^{-53}$ for the corrected-product reflection \eqref{eq:V-reflection}, and $4.21\times10^{-53}$ for the shifted eta equation \eqref{eq:shifted-fe}. At $s=2.37+3.1i$, the residual in the weighted-series Mellin equation \eqref{eq:weighted-Mellin-FE} is approximately $5.01\times10^{-52}$. These are implementation sanity checks of identities already proved symbolically; they are not independent proof and should not be interpreted as evidence for the unproved legacy equation.

\subsection{What the revised slopes do and do not show}

The corrected eta-series calculation shows that a log--log regression of the truncation error at a critical-line zero naturally has slope near $-1/2$. This behavior is predicted by the alternating-tail expansion. It does not establish a functional equation. It is also not the same statistic as the legacy manuscript's regressions of logarithms of the indices where a plotted curve crosses an axis. Since the old raw crossing sequences and their exact finite definitions are unavailable or inconsistent, the old slopes are not reinterpreted here. Any future claim about them must be based on regenerated raw data from a corrected and fully specified algorithm.

\section{The missing correspondence and research questions}\label{sec:open}

The remaining mathematical problem concerns the original variant, not the standard Jacobi kernel and not merely the corrected analytic product $V$.

\begin{question}[Exact theta--eta identification]\label{q:identity}
After a corrected limiting definition of the original variant $W(s)$ has been supplied, is there an explicit meromorphic factor $A(s)$ and a shift $\alpha$ such that
\begin{equation}\label{eq:target-identity}
W(s)=A(s)\LamE(s+\alpha)
\end{equation}
on a nonempty open set? The factor $A$ must be specified independently from the quotient $W(s)/\LamE(s+\alpha)$; otherwise the statement is tautological.
\end{question}

The real-part Euler limit must first be replaced or eliminated, and each finite approximation to $\zeta$ must converge to analytic continuation in the domain being studied. If \eqref{eq:target-identity} is proved with $\alpha=3/8$, the quarter-shift equation for the normalized function follows from Theorem~\ref{thm:shifted-eta}. If the identity is false, the useful outcome would be a counterexample or a different exact/asymptotic description of the constructed sequence.

The next computational study should therefore:
\begin{enumerate}[label=\textbf{Step \arabic*.},leftmargin=3.2em]
  \item freeze a single mathematical definition of the original variant and distinguish its limiting function from each finite approximant;
  \item use a valid complex gamma factor and a convergent or error-controlled representation of $\zeta(s)$;
  \item compare the proposed functional-equation sides at points away from zeros as well as at zeros, recording absolute and relative residuals with precision and truncation depth specified;
  \item repeat calculations with increasing precision and independent evaluation methods;
  \item define the branch of every complex logarithm only after the multiplicative identity is established, avoiding zeros and poles;
  \item release the root list, raw numerical arrays, and plotting script so that regression statistics can be independently reproduced.
\end{enumerate}

\section{Conclusion}
The results have four distinct roles. First, the classical Jacobi kernel has Mellin transform $\Lambda_\eta(s)$ in $\RePart(s)>0$; the kernel splitting also provides its meromorphic continuation, the residue and finite part at $s=0$, several exact moment integrals, and a Mellin inversion formula. A Mellin uniqueness argument proves that the shifted transform $\Lambda_\eta(s+3/8)$ corresponds to the weighted kernel $x^{3/8}K(x)$ within the stated weighted $L^1$ class. Second, the weighted exponential series in the source construction is analytically distinct: its Mellin transform is $B(s)=\Gamma(s)\zeta(s-1)$ and its reflection equation is a meromorphic relation under $s\mapsto3-s$ (center line $\Re s=3/2$), because the zeta factor has shifted argument $s-1$. Its finite sum has an exponentially decaying, explicit tail. Third, the completed eta reflection belongs to a whole family of affine reflections obtained by shifting the argument. The choice $3/8$ yields $1/4-s$, but that observation alone is not diagnostic of the original variant. Fourth, the corrected product is exactly $A(s)\Lambda_\zeta(s)$ with $A(s)=s\Theta^B(s)/2$; it obeys a fully specified $1-s$ functional equation with an elementary meromorphic multiplier. Its Laurent data at $s=0$ and $s=1$, and the cancellation of gamma poles at negative even integers, are also explicit. When both factors are truncated, the product error has the explicit bound in Theorem~\ref{thm:both-truncated}.

The Euler-transformed eta series supplies a second representation with a rigorous geometric truncation bound, although its finite differences require care with cancellation. The stronger alternating-tail expansion explains why ordinary eta-series residuals at a zero of real part $\beta$ have log--log slope approaching $-\beta$. In particular, the slope $-1/2$ on the critical line is predicted by a proved asymptotic, not evidence for a functional equation. The numerical experiments and residual checks are reproducible diagnostics of these stated identities and bounds; they do not constitute independent proofs and do not reconstruct the legacy crossing-index plots.

The original variant remains the unresolved object. Its real-part-based Euler limit does not define a complex function away from the real axis, and its ordinary zeta partial sums do not converge in the critical strip. A new version of the variant must first be defined unambiguously. Only then can one prove or refute an independently specified identity $W(s)=A_*(s)\Lambda_\eta(s+\alpha)$, or derive a different transform law. The auxiliary results here give a mathematically consistent base for that work without claiming that the missing correspondence has already been established.

\appendix
\section{Reproducibility files and numerical conventions}

The supplementary archive contains:
\begin{itemize}[leftmargin=2em]
  \item \texttt{reproduce\_experiments.py}: 50-digit reference evaluation, compensated eta partial sums, error-bound calculations, regressions, and figure generation;
  \item \texttt{eta\_zero\_truncation.csv}: eta partial sums, observed errors, and analytic bounds at $\rho_1$;
  \item \texttt{eta\_euler\_acceleration.csv}: Euler-transformed eta approximants, observed residuals, and geometric bounds at $\rho_1$;
  \item \texttt{fig\_eta\_euler\_acceleration.pdf}: Euler-transformed eta residuals and certified upper bounds;
  \item \texttt{zeta\_eta\_approximation.csv}: eta-based zeta approximations and their propagated bounds at $1/2+82i$;
  \item \texttt{corrected\_product\_error.csv}: error and combined bound for the product with both factors truncated;
  \item \texttt{weighted\_exponential\_truncation.csv}: weighted-series truncation error and its explicit bound;
  \item \texttt{corrected\_product\_zero\_residual.csv}: finite-product residuals at the first nontrivial zeta zero;
  \item \texttt{fig\_eta\_error\_bound.pdf}: eta-tail error and bound;
  \item \texttt{fig\_eta\_scaled\_error.pdf}: scaled eta-tail asymptotic;
  \item \texttt{fig\_product\_zero\_residual.pdf}: product residual at a zero;
  \item \texttt{fig\_zeta\_eta\_error.pdf}: eta-based zeta approximation error;
  \item \texttt{fig\_weighted\_exponential\_error.pdf}: weighted-series truncation;
  \item \texttt{fig\_corrected\_product\_error.pdf}: two-factor product error.
\end{itemize}

Run the experiment from the extracted archive directory with
\begin{lstlisting}[basicstyle=\ttfamily\small,frame=single]
python reproduce_experiments.py
\end{lstlisting}
The script requires Python 3.10 or later and the packages \texttt{mpmath}, \texttt{numpy}, and \texttt{matplotlib}. The regression results are deterministic summaries of the displayed truncation sequence; they are not statistical estimates from random samples.

\end{document}